\documentclass[]{amsart}
\usepackage{times,amsmath}
\usepackage{amssymb}
\usepackage{latexsym}
\newtheorem{theorem}{Theorem}[section]
\newtheorem{lemma}[theorem]{Lemma}
\theoremstyle{definition}
\newtheorem{definition}[theorem]{Definition}

\newtheorem{corollary}[theorem]{Corollary}
\theoremstyle{remark}

\numberwithin{equation}{section}
\newcommand{\As}{\mathcal{A}_{\mathrm{s}}}
\newcommand{\Aw}{\mathcal{A}_{\mathrm{w}}}

\newcommand{\ds}{\mathrm{d}_{\mathrm{s}}}
\newcommand{\dw}{\mathrm{d}_{\mathrm{w}}}

\newcommand{\Dc}{\mathcal{E}}
\newcommand{\Rc}{\widetilde{R}}

\newcommand{\Hw}{H_{\mathrm{w}}}

\newcommand{\ddt}{\frac{d}{dt}}

\newcommand{\ts}{\textstyle}

\begin{document}
\title[Blow-up for the dyadic model]{Blow-up in finite time for the dyadic model of the Navier-Stokes equations.}
\author{Alexey Cheskidov}
\address{University of Michigan, Department of Mathematics, 
Ann Arbor, MI 48109}
\email{acheskid@umich.edu}


\date{January 4, 2006}

\begin{abstract}
We study the dyadic model of the Navier-Stokes equations
introduced by Katz and Pavlovi\'c. They showed a finite time
blow-up in the case where the dissipation degree $\alpha$ is less than $1/4$.
In this paper we prove the existence of weak solutions for all $\alpha$,
energy inequality for every weak solution with nonnegative initial data
starting from any time, 
local regularity for $\alpha > 1/3$, and global regularity for $\alpha \geq 1/2$.
In addition, we prove a finite time blow-up in the case where $\alpha<1/3$.
It is remarkable that the model with $\alpha=1/3$ enjoys the same estimates
on the nonlinear term as  the 4D Navier-Stokes equations.
Finally, we discuss a weak global attractor,
which coincides with a maximal bounded invariant set for all $\alpha$ and
becomes a strong global attractor for $\alpha \geq 1/2$.
\end{abstract}

\maketitle

\section{Introduction}
The regularity of the 3D incompressible Navier-Stokes equations (NSE) remains
a significant problem. This, among many other
open problems connected with the 3D NSE, depends on the estimates on the
inertial term $(\mathbf{u} \cdot \nabla)\mathbf{u}$ in the equations. In this paper
we study a dyadic model, which has similar properties
to the 3D NSE,  the same estimates on the inertial term,
and the same open question concerning the regularity of the solutions.

There have been many simple models
proposed in the literature that capture some  essential features of
the 3D NSE. Among these are shell models
of turbulence, which have been investigated for many years
(see \cite{B,F,G,LPPPV,OY}).
Recently, some of these models, as well as some new ones, were extensively
studied analytically. In \cite{CLT}, Constantin, Levant, and Titi study the "sabra"
shell model of turbulence, proving a global regularity and the existence of a
finite dimensional global attractor and inertial manifold.

In \cite{KP}, Katz and Pavlovi\'c introduced another shell-type model,
the dyadic model for the Euler and
Navier-Stokes equations. This model, motivated by \cite{KP1}, is an infinite
system of nonlinear
ODEs that describes evolutions of wavelet coefficients. In \cite{FP},
Friedlander and Pavlovi\'c proposed a three-dimensional vector model for
the Euler equations, similar to a quasi-linear approximation of the 3D
Navier-Stokes system constructed by Dinaburg and Sinai \cite{DS}. Both of these
dyadic models can be reduced to the following system of nonlinear ODEs:
\begin{equation} \label{e:dintro}
\ddt u_n + \nu \lambda^{2\alpha n}u_n - \lambda^n u_{n-1}^2 + \lambda^{n+1} u_n u_{n+1}
=g_n, \qquad n\in \mathbb{N},
\end{equation}
where $u_0=0$. Here, $\lambda >1$, $\nu\geq 0$ is the viscosity, and $\alpha >0$ is
the dissipation degree. Note that we also include a force $g=(g_1,g_2,\dots)$
in the model.
For $u=(u_1,u_2,\dots)$, the dyadic model can be written as
\begin{equation*}
\ddt u + \nu Au + B(u,u) = g,\\
\end{equation*}
where
\[
(Au)_n = \lambda^{2\alpha n} u_n, \qquad
(B(u,u))_n=
- \lambda^n u^2_{n-1} + \lambda^{n+1} u_n u_{n+1},
\]
and $u_0=0$. 
Note that for $\alpha=2/5$ the following are 
sharp estimates on the inertial term:
\[
|(B(u,u),Au)| \lesssim  |Au|^{3/2}|A^{1/2}u|^{3/2},
\]
where $(\cdot,\cdot)$ and $|\cdot|$ are the $l^2$-inner product and norm
respectively. The best known estimates on the inertial term of the 3D NSE
are the same, with $(\cdot,\cdot)$ and $|\cdot|$ being the $L^2$-inner product
and norm.

Katz and Pavlovi\'c \cite{KP}
proved that under certain assumptions on
the initial conditions,  solutions of the inviscid dyadic model blow up in
finite time in a norm stronger than the $l^2$-norm.
Later, Waleffe \cite{W} derived the inviscid dyadic model from the Burgers
equation reducing it to \eqref{e:dintro} with $\nu =0$.
Recently, Kiselev and Zlato\v s \cite{KZ} sharpened the blow-up result
for the dyadic model and studied a very similar Obuhov model
(see \cite{O}) proving a global regularity of every solution with regular initial data.
In addition, the existence of a global attractor of the
inviscid dyadic model is proved in \cite{CFP1}. This surprising fact is a result of a self-dissipation
mechanism due to the loss of regularity of solutions.

In this paper we will study the viscous dyadic model, i.e., the model~\eqref{e:dintro} with $\nu>0$.
In this case Katz and Pavlovi\'c \cite{KP} obtained a finite time blow-up
of  solutions with certain initial data  when $\alpha <1/4$.
Our main goal is to prove a finite time blow-up in the case where $\alpha<1/3$.
It is remarkable that in the critical case $\alpha =1/3$ the following are
sharp estimates on the nonlinear term:
\[
|(B(u,u),Au)| \lesssim  |Au|^{2}|A^{1/2}u|.
\]
Note that the best known estimates
on the nonlinear term of the 4D NSE are exactly the same.

This paper is structured as follows. We start with surveying the dyadic model
in Section~\ref{s:dyadic}. In Section~\ref{s:setting} we introduce a functional
setting and define weak and strong solutions. In Section~\ref{s:apriori}
we derive some {\it a priori} estimates and prove the
existence of weak solutions to the dyadic model.
This is done by taking a limit of the Galerkin approximation, which also
results in the energy inequality for a limit solution (which might not
be unique) for almost all time. Then we show that the lack of backward
energy transfer implies that every weak solution with nonnegative initial
data satisfies the energy
inequality starting from any time. Finally, using a classical NSE technique,
we show a local regularity for $\alpha >1/3$ and a global regularity for
$\alpha \geq 1/2$.

In Section~\ref{s:blowup}, inverting the Sobolev-type estimates, we prove that
every solution with large $H^{\epsilon}$-norm blows up in finite time in
$H^{1/3+\epsilon}$-norm, $\epsilon>0$. Here, $H^\gamma$-norm is defined as
$\|u\|_\gamma=(\sum \lambda^{2\gamma n}u_n^2)^{1/2}$.
Note that we also use such a technique in \cite{CF},
where we study a similar model that has coefficients growing as power
functions. That model is introduced as an example of a NSE-like
dynamical system that possesses a weak global attractor, on which
all the solutions blow up in finite time. The reason for the power-law
growth for the coefficients is to mimic the growth of eigenvalues of
the Stokes operator in 3D. It is remarkable that due to a slower growth of
the coefficients, the model in \cite{CF} still possesses a gap between the regions
of a local regularity and finite time blow-up.

Lastly, in Section~\ref{s:attractor} we discuss a weak global attractor
for the dyadic model. The weak
global attractor is the minimal weakly closed weakly attracting set. 
Using results from \cite{C,CF}, we show that the weak global attractor is also the
maximal bounded invariant set. Moreover, for $\alpha \geq 1/2$ all the
trajectories are continuous in $l^2$, which implies that the weak global
attractor is in fact a strongly compact strong global attractor.

Note that there is still a gap between the regions of global regularity and
blow-up in finite time, which means that the developed technique is not sharp
enough to separate these two behaviors. 
Since most of the proofs in the theory of the Navier-Stokes equations go through for the dyadic model, 
a better understanding of the dyadic or similar shell models might provide
insight into the regularity problem for the Navier-Stokes equations.

\section{Dyadic model} \label{s:dyadic}

Here we recall a derivation of the dyadic model for the equations of fluid
motion by Katz and Pavlovi\'c \cite{KP}. A cube $Q \subset \mathbb{R}^3$ is called
dyadic if its side length is $2^l$, and the corners are on the
lattice $2^l\mathbb{Z}^3$, for some integer $l$. For a dyadic cube
$Q$ with side length $2^{-j}$, its parent $\tilde{Q}$ is a unique
dyadic cube with side length $2^{-j+1}$ that contains $Q$. For
$m\geq 1$, let $C^m(Q)$ be the set of all $m$th order grandchildren
of $Q$, i.e., all the dyadic cubes with side length $2^{-j-m}$
that are contained in $Q$. For instance, $C^1(Q)$ consists of
$2^3$ children of $Q$.

Now a scalar-valued function $u(x,t)$ can be represented by the following
wavelet expansion
\[
u(x,t) = \sum_Q u_Q(t) w_Q(x),
\]
where $\{w_Q\}$ is an orthonormal in $L^2(\mathbb{R})$ family of
wavelets, such that $w_Q$ is localized on $Q$.
Define the Laplacian  in the following way:
\[
\Delta u=\sum_Q 2^{2j(Q)}u_Q(t)\omega_Q(x),
\]
where $2^{-j(Q)}$ is a side length of a dyadic cube $Q$.
Katz and Pavlovi\'c define the cascade operator as follows:
\[
(C(u,v))_Q=-2^\frac{5j}{2}u_{\tilde{Q}} v_{\tilde{Q}}+
2^\frac{5(j+1)}{2}u_Q\sum_{Q'\in C^1(Q)} v_{Q'},
\]
where $Q$ is a dyadic cube with side length $2^{-j}$.
The dyadic Navier-Stokes equation with hypo-dissipation is written
as
\[
\ddt u +C(u,u) + \nu (\Delta)^\alpha u =0,
\]
where we include viscosity $\nu$, which is chosen to be one in \cite{KP}.
In terms of the wavelet coefficients $u_Q$, this equation can be
written as
\[
\ddt u_Q(t) = -\nu 2^{2\alpha j}u_Q(t)
+2^{\frac{5j}{2}}u_{\tilde{Q}}^2(t)
-2^{\frac{5(j+1)}{2}}u_Q \sum_{Q'\in C^1(Q)} u_{Q'}(t).
\]

As it was proposed in \cite{W} in the case $\nu=0$, we simplify the model
in the following way.
Let $Q_1$ be a dyadic cube with side length $2^{-1}$. Let $v_1(t) =u_{Q_1}(t)$
and $v_{j}(t)=u_{Q_j}(t)$, where $Q_j$ is some dyadic cube
in $C^{j-1}(Q_1)$, $j\geq 2$.
We will only consider the initial conditions for which
$u_Q(t_0)= v_{j}(t_0)$ for all cubes $Q \in C^{j-1}(Q_1)$ for $j\geq 2$, and
$u_Q(t_0)=0$ for all dyadic cubes with side length larger than $2^{-1}$.
Then, for every $j\geq 2$, we have $u_Q(t)=v_{j}(t)$ for all $Q \in C^{j-1}(Q_1)$.
Now, denoting $v_0=0$, we
obtain the following system of equations for $v_j(t)$:
\[
\ddt v_j(t) = -\nu 2^{2\alpha j}v_j(t) +2^{\frac{5j}{2}}v_{j-1}^2 - 2^32^{\frac{5(j+1)}{2}} v_j v_{j+1} , \qquad j \geq 1,
\]
Finally, the change of variables
\[
u_j(t)=2^{\frac{3j}{2}} v_j(t/8).
\]
reduces the equations to
\[
\ddt u_j = -\tilde{\nu} 2^{2\alpha j} u_j + 2^j u_{j-1}^2 - 2^{j+1}u_j u_{j+1},
\qquad j\geq 1,
\]
with $\tilde{\nu} =\nu/8$.

\section{Functional setting} \label{s:setting}
Let us denote $H=l^2$ with the usual scalar product and norm:
\[
(u,v):= \sum_{n=1}^{\infty} u_nv_n, \qquad |u|:=\sqrt{(u,u)}.
\]
The norm $|u|$ will be called the energy norm.
Let $A:D(A) \to H$ be the Laplace operator defined by
\[
(Au)_n = \lambda^{2\alpha n} u_n, \qquad n\geq 1,
\]
for some $\lambda >1$.
The domain $D(A)$ of this operator is a dense subset of $H$.
Note that $A$ is a positive definite operator
whose eigenvalues are
\[
0<\lambda^{2\alpha} \leq \lambda^{4\alpha} \leq \lambda^{6\alpha} \leq \dots
\]
Let $H^{\gamma}=A^{-\gamma/(2\alpha)}H$ endowed with the
following scalar product and norm: 
\[
((u,v))_{\gamma}:= \sum_{n=1}^{\infty} \lambda^{2\gamma n}u_nv_n, \qquad \|u\|_{\gamma}:=\sqrt{((u,u))_\gamma}.
\]
In the special case $\gamma=\alpha$, let $V=H^{\alpha} =A^{-1/2}H$ and
\[
((u,v)):=((u,v))_{\alpha}, \qquad \|u\|:=\|u\|_\alpha.
\]
This double norm $\|u\|$ will be called the enstrophy norm.
Note that we have an equivalent of the Poincar\'{e} inequality
\[
|u|^2 \leq \frac{1}{\lambda^{2\alpha}} \|u\|^2.
\]
Let also
\[
\ds(u,v):=|u-v|, \qquad
\dw(u,v):= \sum_{n=1}^\infty \frac{1}{2^{(n^2)}}
\frac{|u_n-v_n|}{1 + |u_n-v_n|},
\qquad u,v \in H.
\]
Here, $\ds$ is a strong distance, and $\dw$ is a weak distance
that induces a weak topology on any bounded subset of $H$. Hence,
a bounded sequence $\{u^k\} \subset H$ converges to $u\in H$ weakly,
i.e.,
\[
\lim_{k\to \infty} (u^k, v) = (u,v),  \qquad \forall v \in H,
\]
if and only if
\[
\dw(u^k, u) \to 0 \qquad \text{as} \qquad  k\to \infty.
\]
We also recall that if $u^k \to u$ weakly in $H$ as $k\to \infty$, then
\[
\liminf_{k \to \infty}|u^k| \geq |u|.
\]
Let
\[
C([0,T];\Hw):=\{u(\cdot): [0,T] \to H, u_n(t) \mbox{ is continuous for all }
n\}
\]
endowed with the distance
\[
d_{C([0, T];\Hw)}(u,v) = \sup_{t\in[0,T]}\dw(u(t),v(t)). 
\]
Let also
\[
C([0,\infty);\Hw):=\{u(\cdot): [0,\infty) \to H, u_n(t) \mbox{ is continuous for all }
n\}
\]
endowed with the distance
\[
d_{C([0, \infty);\Hw)}(u,v) = \sum_{T\in \mathbb{N}} \frac{1}{2^T} \frac{\sup\{\dw(u(t),v(t)):0\leq t\leq T\}}
{1+\sup\{\dw(u(t),v(t)):0\leq t\leq T\}}.
\]

In this paper, the dyadic model of the Navier-Stokes equations will
be written as
\begin{equation} \label{model}
\left\{
\begin{aligned}
&\ddt u_n + \nu \lambda^{2\alpha n}u_n - \lambda^n u_{n-1}^2 + \lambda^{n+1} u_n u_{n+1}
=g_n, \qquad n=1,2,3\dots\\
&u_0=0,
\end{aligned}
\right.
\end{equation}
for some parameter $\lambda >1$, the viscosity $\nu >0$, the
dissipation degree $\alpha>0$, and the force $g=(g_1,g_2,\dots)$.
For simplicity, assume that $g$ is independent of time, $g\in H$, and
$g_n\geq 0$ for all $n$.

For $u=(u_1,u_2,\dots)$, the dyadic model can be written in a more condensed
form as
\begin{equation} \label{model1}
\ddt u + \nu Au + B(u,u) = g,\\
\end{equation}
where
\[
(B(u,v))_n=
\left\{
\begin{array}{ll}
- \lambda^n u_{n-1} v_{n-1} + \lambda^{n+1} u_n v_{n+1}, & n=2,3,\dots\\
\lambda^{2} u_1 v_{2}, & n=1.
\end{array}
\right.
\]
Clearly, the bilinear operator $B$ enjoys the orthogonality property:
\[
\begin{split}
\left(B(u,v),v\right)&=\sum_{n=1}^{\infty}\left(
-\lambda^nu_{n-1} v_{n-1} v_n + \lambda^{n+1} u_n v_{n+1} v_n \right)\\
&=
\sum_{n=1}^{\infty}\left(
-\lambda^nu_{n-1} v_{n-1} v_n +\lambda^n u_{n-1} v_n v_{n-1} \right)\\
&=0.
\end{split}
\]
Note that we always use a convention that $u_0=0$.
\begin{definition}
A weak solution on $[T,\infty)$ (or $(-\infty, \infty)$, if
$T=-\infty$) of \eqref{model} is an $H$-valued
function $u(t)$ defined for $t \in [T, \infty)$, such that $u_n \in C^1([T,\infty))$
and $u_n(t)$ satisfies \eqref{model} for all $n$.
\end{definition}

Note that since $(B(u,u))_n$ has a finite number of terms, the notions of
a weak solution and a classical solution (of a system of ODEs) coincide.
Hence, the weak solutions will be often
called solutions in the remainder of the paper. 
Note that if $u(t)$ is a solution on $[T,\infty)$, then automatically
$u_n \in C^{\infty}([T,\infty))$. We say that a solution $u(t)$ is
strong (or regular) on some interval $[T_1,T_2]$, if $\|u(t)\|$ is bounded on $[T_1,T_2]$.
A solution is strong on $[T_1,\infty)$, if it is strong on every interval $[T_1,T_2]$, $T_2\geq T_1$.

\begin{definition} \label{d:ex}
A Leray-Hopf solution of \eqref{model} on the interval $[T, \infty)$
is a weak solution of \eqref{model} on $[T,\infty)$ satisfying the
energy inequality
\begin{equation} \label{EI}
|u(t)|^2 + 2\nu \int_{t_0}^t \|u(\tau)\|^2 \, d\tau \leq
|u(t_0)|^2 + 2\int_{t_0}^t (g, u(\tau)) \, d\tau
\end{equation}
 for all $T \leq t_0 \leq t$,
$t_0$ a.e. in $[T,\infty)$. The set $Ex$ on which the energy
inequality does not hold will be called the exceptional set.
\end{definition}

Note that
the complement of the exceptional set $Ex$ coincides with the set of points of
strong continuity from the right.
Later we will prove that every solution $u(t)$ with $u_n(T)\geq0$
is a Leray-Hopf solution on $[T,\infty)$, and that the energy inequality for such
a solution is satisfied starting from any time $t_0 \geq T$, i.e., $Ex = \emptyset$.

\section{A priori estimates and existence of weak and strong solutions}
\label{s:apriori}
We start with some {\it a priori} estimates. 
\\

\noindent
{\bf Energy estimates.} Formally taking a scalar product of the equation~\eqref{model} with $u$, we obtain
\[
\begin{split}
\frac{1}{2} \ddt |u|^2  &\leq  -\nu \|u\|^2 + |g||u|\\
&\leq  -\nu |u |^2 + \frac{\nu}{2}|u|^2 + \frac{|g|^2}{2\nu}\\
&= -\frac{\nu}{2}|u|^2 + \frac{|g|^2}{2\nu}.
\end{split}
\]
Using Gronwall's inequality, we conclude that
\begin{equation} \label{eq:defK}
|u(t)|^2 \leq e^{- \nu t} |u(0)|^2 +
\frac{|g|^2}{\nu^2}(1-e^{-\nu t}).
\end{equation}
Hence, $B=\{u\in H: \ |u| \leq R\}$ is an absorbing ball for the
Leray-Hopf solutions, where $R$ is any number larger that $|g|/\nu$.
Note that this result will later follow rigorously  from the energy inequality.

Next, taking a limit of the Galerkin approximation, we will
prove the existence of Leray-Hopf solutions to \eqref{model}.

\begin{theorem} \label{thm:Leray}
For every $u^0 \in H$ and $g \in H$,
there exists a solution of \eqref{model} with $u(0)=u^0$.
Moreover, the energy inequality
\begin{equation*} 
|u(t)|^2 + 2\nu \int_{t_0}^t \|u(\tau)\|^2 \, d\tau \leq
|u(t_0)|^2 + 2\int_{t_0}^t (g, u(\tau)) \, d\tau
\end{equation*}
holds for all $0 \leq t_0 \leq t$, $t_0$ a.e. in $[0,\infty)$.
\end{theorem}
\begin{proof}
Let $u^0 \in H$.
We will show the existence of a weak solution by taking a limit of the Galerkin
approximation $u^k(t)=(u^k_1(t),\dots,u^k_k(t),0,0,\dots)$ with
$u^k_n(0)=u^0_n$ for $n=1,2,...,k$, which satisfies
\begin{equation} \label{galerkin}
\left\{
\begin{aligned}
&\ddt u^k_n + \nu \lambda^{2\alpha n}u^k_n - \lambda^n (u^k_{n-1})^2 +
\lambda^{n+1} u^k_n u^k_{n+1}
=g_n, \qquad n\leq k-1,\\
&\ddt u^k_k + \nu \lambda^{2\alpha k}u^k_k - \lambda^k (u^k_{k-1})^2
=g_k,
\end{aligned}
\right.
\end{equation}
where $u^k_0=0$.
First, note that the energy estimate~\eqref{eq:defK} obviously holds for $u^k(t)$.
Hence, from the theory of ordinary differential equations we know that 
there exists a unique solution $u^k(t)$ to \eqref{galerkin} on $[0,\infty)$.
Next, we will show that a sequence of the Galerkin approximations $\{u^k\}$ is
weakly equicontinuous. Indeed, thanks to the energy estimate~\eqref{eq:defK},
there exists $M$, such that
\[
u^k_n(t) \leq M, \qquad \forall n, k, t\geq 0.
\]
Therefore,
\[
\begin{split}
|u^k_n(t) - u^k_n(s)| &\leq \left| \int_s^t \left(-\nu \lambda^{2\alpha n} u_n^k
+\lambda^n(u^k_{n-1})^2 -\lambda^{n+1}u^k_n u^k_{n+1} + g_n\right) \, d\tau\right|\\
&\leq (\nu\lambda^{2\alpha n}M +\lambda^n M^2 + \lambda^{n+1}M^2 + g_n)
|t-s|,
\end{split}
\]
for all $n$, $k$, $t\geq 0$, $s\geq0$. Thus,
\[
\begin{split}
\dw(u^k(t), u^k(s))&= \sum_{n=1}^\infty \frac{1}{2^{(n^2)}}
\frac{|u^k_n(t)-u^k_n(s)|}{1 + |u^k_n(t)-u^k_n(s)|}\\
&\leq c|t-s|,
\end{split}
\]
for some constant $c$ independent of $k$. Hence, $\{u^k\}$ is an equicontinuous sequence
of functions in $C([0,\infty);\Hw)$ with bounded initial data. Therefore,
the Ascoli-Arzela theorem implies that
$\{u^k\}$ is relatively compact in $C([0,T]; \Hw)$ for every $T\geq0$. By
a diagonalization process it follows that $\{u^k\}$ is relatively compact in
$C([0,\infty); \Hw)$. Hence, passing to a subsequence,
we obtain that there exists a weakly continuous $H$-valued
function $u(t)$, such that
\begin{equation} \label{e:weakconv}
u^{k_j} \to u \qquad \mbox{as} \qquad k_j \to \infty  \qquad \mbox{in}
\qquad  C([0,\infty); \Hw).
\end{equation}
In particular, $u^{k_j}_n(t) \to u_n(t)$ as $k_j \to \infty$, for all $n, t \geq 0$.
Thus, $u(0) = u^0$.
In addition, note that
\[
u^{k_j}_n(t) = u^{k_j}_n(0) + \int_0^t (-\nu \lambda^{2\alpha n}u^{k_j}_n +
\lambda^n (u^{k_j}_{n-1})^2 - \lambda^{n+1} u^{k_j}_n u^{k_j}_{n+1} +g_n) \, d\tau,
\]
for $n \leq k_j-1$. Taking the limit as $k_j \to \infty$, we obtain
\[
u_n(t) = u_n(0) + \int_0^t (-\nu \lambda^{2\alpha n}u_n +
\lambda^n u_{n-1}^2 - \lambda^{n+1} u_n u_{n+1} +g_n) \, d\tau.
\]
Since $u_n(t)$ is continuous, it follows that
$u_n \in C^1([0,\infty))$ and satisfies
\eqref{model}.

It remains to prove that $u(t)$ satisfies the energy inequality.
Note that $u^{k_j}(t)$ satisfies the energy equality
\[
|u^{k_j}(t)|^2 + 2\nu \int_{t_0}^t \|u^{k_j}(\tau)\|^2 \, d\tau =
|u^{k_j}(t_0)|^2 + 2\int_{t_0}^t (g, u^{k_j}(\tau)) \, d\tau,
\]
for all $t \geq t_0\geq0$. Hence, the sequence $\{u^{k_j}\}$ is bounded
in $L^2([t_0,t];V)$ for all $t \geq t_0\geq0$. 
This, together with \eqref{e:weakconv}, implies that
\[
\int_{t_0}^t |u^{k_j}(\tau)-u(\tau)|^2 \, d\tau \to 0, \qquad \mbox{as}
\qquad  k_j \to \infty,
\]
for all $t\geq t_0\geq 0$.
In particular, $|u^{k_j}(t)| \to |u(t)|$ as $k_j \to \infty$ a.e. in $[0,\infty)$.
Take any $t_0\geq 0$ for which $|u^{k_j}(t_0)| \to |u(t_0)|$ as $k_j \to \infty$.
For every $N\geq 0$, we have
\[
|u^{k_j}(t)|^2 + 2\nu \int_{t_0}^t \sum_{n\leq N} \lambda^{2\alpha n} u^{k_j}_n
(\tau)^2 \, d\tau \leq
|u^{k_j}(t_0)|^2 + 2\int_{t_0}^t (g, u^{k_j}(\tau)) \, d\tau.
\]
Since $u^{k_j}(t) \to u(t)$ weakly in $H$ as $k_j \to \infty$
for all time $t\geq 0$, we have that
\[
|u(t)|^2 + 2\nu \int_{t_0}^t \sum_{n\leq N} \lambda^{2\alpha n} u_n(\tau)^2 \, d\tau \leq
|u(t_0)|^2 + 2\int_{t_0}^t (g, u(\tau)) \, d\tau.
\]
Finally, taking the limit as $N \to \infty$ and using Levi's
theorem, we obtain
\[
|u(t)|^2 + 2\nu \int_{t_0}^t \|u(\tau)\|^2 \, d\tau \leq
|u(t_0)|^2 + 2\int_{t_0}^t (g, u(\tau)) \, d\tau,
\]
for all $0 \leq t_0 \leq t$, $t_0$ a.e. in $[0,\infty)$.
\end{proof}

Note that this was a classical proof from the theory of the NSE.
Using the fact that there is no backward energy transfer, we can actually
show that every solution with $u_n(0)\geq 0$ is a Leray-Hopf solution
and, moreover, is continuous from the right in $H$ for all time. 

\begin{theorem} \label{t:iquality}
Let $u(t)$ be a solution of \eqref{model} with $u_n(0)\geq 0$.
Then $u_n(t) \geq 0$ for all $t>0$, and 
$u(t)$ satisfies the energy inequality
\begin{equation} \label{ee}
|u(t)|^2 + 2\nu \int_{t_0}^t \|u(\tau)\|^2 \, d\tau \leq
|u(t_0)|^2 + 2\int_{t_0}^t (g, u(\tau)) \, d\tau,
\end{equation}
for all $0 \leq t_0 \leq t$.
\end{theorem}
\begin{proof}
A general solution for $u_n(t)$ can be written as
\begin{multline} \label{e:gensol}
u_n(t)=u_n(0)\exp\left(-
\int_{0}^t\nu \lambda^{2\alpha n} + \lambda^{n+1} u_{n+1}(\tau) \, d\tau\right)\\
+ \int_{0}^t\exp\left(-\int_{s}^{t}\nu \lambda^{2\alpha n}+\lambda^{n+1}
u_{n+1}(\tau) \, d\tau\right) (g_n+\lambda^n u_{n-1}^2(s) )\, ds.
\end{multline}
Recall that $g_n \geq 0$ for all $n$.
Since $u_n(0) \geq 0$ for all $n$, then  $u_n(t)\geq 0$ for all $n$, $t>0$.
Hence, multiplying \eqref{model} by $u_n$,
taking a sum from $1$ to $N$, and integrating between $t_0$ and $t$,
we obtain
\[
\begin{split}
\sum_{n=1}^N u_n(t)^2  -& \sum_{n=1}^N u_n(t_0)^2
 +2\nu \int_{t_0}^t \sum_{n=1}^N \lambda^{2\alpha n}u_n(\tau)^2
 \, d\tau \\
&= - 2\int_{t_0}^t  \lambda^{N+1} u_N^2
u_{N+1} \, d\tau  +2\int_{t_0}^t\sum_{n=1}^N  g_n u_n \, d\tau\\
&\leq 2\int_{t_0}^t\sum_{n=1}^N g_n u_n \, d\tau.
\end{split}
\]
Taking the limit as $N \to \infty$, we obtain \eqref{ee}.

\end{proof}

\noindent
{\bf Enstrophy estimates.}
We obtain the following estimate for the nonlinear term:
\[
\begin{split}
| (B(u,u),Au) | &=
\left| \sum_{n=1}^\infty \left[\lambda^{2\alpha (n+1)}-\lambda^{2\alpha n} \right]
\lambda^nu_n^2 u_{n+1} \right|\\
&=
(\lambda^{\alpha} - \lambda^{-\alpha})
\left| \sum_{n=1}^\infty \lambda^{(\alpha+1) n} u_n^2 \lambda^{\alpha(n+1)} u_{n+1}
\right|\\
&\leq c_\mathrm{b} (\max_n |\lambda^{\alpha n}u_n|) \sum_{n=1}^{\infty}
\lambda^{(\alpha+1) n} u_n^2\\
&\leq c_\mathrm{b} \|u\| \sum_{n=1}^{\infty} \lambda^{(\alpha+1)n} u_n^2,
\end{split}
\]
where $c_\mathrm{b} =\lambda^\alpha -\lambda^{-\alpha} >0$.
When $\alpha \in[1/3, 1]$, H\"older's inequality implies
\[
\begin{split}
| (B(u,u),Au) | &\leq c_\mathrm{b}\|u\| |Au|^{1/\alpha -1} |A^{1/2}u|^{3-1/\alpha}\\
&= c_\mathrm{b}|Au |^{1/\alpha-1}
\|u\|^{4-1/\alpha}.
\end{split}
\]
Choosing $u$ to have only two consecutive nonzero terms, it is easy to
check that these estimates are sharp. Moreover, when $\alpha=2/5$, we have
\[
| (B(u,u),Au) | \leq c_\mathrm{b} |Au |^{3/2}\|u\|^{3/2},
\]
which is the same as the Sobolev estimate for the inertial term of the 3D NSE
(see, e.g., \cite{ConF,T1}).
Therefore, taking a scalar product of the equation \eqref{model} with $Au$
and using Young's inequality, we obtain
\[
\begin{split}
\frac{1}{2}\ddt \|u\|^2 &\leq -\nu|Au|^2 + c_\mathrm{b}|Au |^{3/2}
\|u\|^{3/2} + (g, Au)\\
&\leq -\nu|Au|^2 + \frac{\nu}{3} |Au|^2 + \frac{3^6c_\mathrm{b}^4}{2^8\nu^3} \|u\|^6 +
\frac{3}{4\nu} |g|^2 + \frac{\nu}{3} |Au|^2\\
&= -\frac{\nu}{3}|Au|^2 + \frac{3^6c_\mathrm{b}^4}{2^8 \nu^3} \|u\|^6 +
\frac{3}{4\nu} |g|^2,
\end{split}
\]
a Riccati-type inequality for $\|u\|^2$. Hence, the model has the same
enstrophy estimate as the 3D NSE, similar properties, and the same open question concerning the regularity of the solutions in the case
$\alpha = 2/5$.

Another interesting case is $\alpha=1/3$.  Then we have
\begin{equation} \label{e:est}
| (B(u,u),Au) | \leq c_\mathrm{b}|Au |^{2}\|u\|,
\end{equation}
which corresponds to the 4D Navier-Stokes equations.

Now consider the case where $\alpha > 1/3$. Formally, we have
\[
\begin{split}
\frac{1}{2}\ddt \|u\|^2 &\leq -\nu|Au|^2 +  c_\mathrm{b}|Au |^{1/\alpha-1}\|u\|^{4-1/\alpha}
+ (g,Au)\\
&\leq -\nu|Au|^2 +  \frac{\nu}{3}|Au |^2 + c \|u\|^{\frac{8\alpha-2}{3\alpha-1}} + \frac{3}{4\nu}|g|^2 +\frac{\nu}{3}|Au|^2\\
& = - \frac{\nu}{3} |Au|^2 + c\|u\|^{\frac{8\alpha-2}{3\alpha-1}} + \frac{3}{4\nu}|g|^2,
\end{split}
\]
for some constant $c >0$. This means that if  the initial data is in $V$,
then $u(t)$ remains bounded in $V$ for some time $T$.
Applying the above estimate to the Galerkin approximation and taking a limit,
we immediately obtain the following local regularity result.
\begin{theorem}
If $\alpha >1/3$, then for any $u^0 \in V$ there exists a strong solution $u(t)$ to \eqref{model}
on some time interval $[0,T]$, $T>0$ with $u(0)=u^0$. 
\end{theorem}

Finally, consider the case $\alpha \geq 1/2$. In this case the enstrophy
estimate implies
\[
(B(u,u),Au) \leq c_\mathrm{b}|Au| \|u\|^2.
\]
Therefore, formally, we have
\[
\begin{split}
\frac{1}{2}\ddt \|u\|^2 &\leq -\nu|Au|^2 +  c_\mathrm{b}|Au| \|u\|^2 + (g,Au)\\
&\leq -\nu|Au|^2 +  \frac{\nu}{3}|Au |^2 + \frac{3c_\mathrm{b}^2}{4\nu} \|u\|^4 + \frac{3}{4\nu}|g|^2 +\frac{\nu}{3}|Au|^2\\
& = - \frac{\nu}{3} |Au|^2 + \frac{3c_\mathrm{b}^2}{4\nu} \|u\|^4 + \frac{3}{4\nu}|g|^2.
\end{split}
\]
This is again a Ricatti-type inequality. Assume that $u(t)$ is a strong solution
on some interval $(0, t^*)$, and $\|u(t)\| \to +\infty$ as $t \to t^*-$. Then
\[
\|u(t)\|^2 \geq \frac{c}{t^*-t}, \qquad 0<t<t^*,
\]
for some positive constant $c$. However, this means that $\|u(t)\|^2$ is
not locally integrable, which is in contradiction with the energy inequality.
Hence, if the initial data $u^0 \in V$, then $\|u(t)\|$ is bounded
on every interval $[0, T]$, $T>0$, and we have the following.
\begin{theorem} \label{t:gregularity}
If $\alpha \geq 1/2$, then for any $u^0 \in V$ there exists a strong solution $u(t)$ to \eqref{model}
on $[0,\infty)$ with $u(0)=u^0$. 
\end{theorem}

\section{Blow-up in finite time} \label{s:blowup}

Let $\alpha <1/3$ and $\gamma >0$. In this section we will prove
that every solution $u(t)$ with large enough $\|u(0)\|_\gamma$ blows up
in finite time in $H^{1/3+\gamma}$ norm.
The idea is the following. Taking a scalar product of the equation with
$A^{\gamma/\alpha} u$, we obtain
\[
\frac{1}{2}\ddt \|u\|_\gamma^2 = -\nu\|u\|^2_{\alpha+\gamma}+(B(u,u),A^{\gamma/\alpha} u) +
 (g,A^{\gamma/\alpha} u).
\]
In order to show a blow-up, we, in some sense, will invert the Sobolev
estimates for the nonlinear term. Note that
\[
(B(u,u),A^{\gamma/\alpha} u) \sim \sum_n \lambda^{(1+2\gamma)n} u_n^2 u_{n+1}.
\]
If $u_n\geq0$ is monotonically decreasing in $n$, then
\[
(B(u,u),A^{\gamma/\alpha} u) \gtrsim \sum_n
 \lambda^{(1+2\gamma)n} u_n^3.
\]
Obviously, this is not true in general. For example, if $u_n=0$ for
even $n$, then $(B(u,u),A^{\gamma/\alpha} u)=0$. However, we will prove that
a similar estimate holds if we use the following function instead of
$H^\gamma$-norm:
\[
H(t) :=  \|u(t)\|_\gamma^2+ c\sum_n \lambda^{ 2\gamma n}(u_nu_{n+1})(t),
\]
for some constanat $c>0$.
More precisely, we will show that if $\|u\|_\gamma$ is large enough, then

\[
\begin{split}
\frac{1}{2}\ddt H  &\gtrsim -\nu\|u\|^2_{\alpha +\gamma} +
\sum_n \lambda^{(1+2\gamma)n}u_n^3\\
&\gtrsim -\nu\|u\|^2_{\alpha+\gamma} + \|u\|_{\alpha+\gamma} ^3\\
&\gtrsim H^{3/2},
\end{split}
\]
provided that $\alpha < 1/3$. We will start with the following estimate.

\begin{lemma} \label{l:jensen}
If $\alpha<1/3$, then for any $\gamma \in(0,1-3\alpha)$ there
exists a positive constant $A$, such that  
\[
\sum_{n=1}^\infty \lambda^{(1+2\gamma)n}|u_n|^3 \geq
A \|u\|_{\alpha + \gamma}^3.
\]
\end{lemma}
\begin{proof}
Let $\epsilon:=2-6\alpha - 2\gamma>0$. Note that
$\lambda^{-\epsilon} < 1$.
Let
\[
A(\gamma) := \left( \sum_{n=1}^\infty \lambda^{-\epsilon n}\right)^{-1/2} =
\sqrt{\lambda^{\epsilon}-1}.
\]
H\"older's inequality with $p=3$ and $q=3/2$ implies
\[
\begin{split}
\|u\|_{\alpha + \gamma}^2
&= \sum_{n=1}^\infty \lambda^{2(\alpha +\gamma)n} u_n^2\\
&\leq \left(  \sum_{n=1}^\infty \lambda^{-\epsilon n}  \right)^{1/3}
 \left(  \sum_{n=1}^\infty \lambda^{(1+2\gamma) n}|u_n|^3 \right)^{2/3}\\
&= A^{-2/3}
\left(  \sum_{n=1}^\infty \lambda^{(1+2\gamma) n}|u_n|^3 \right)^{2/3}.
\end{split}
\]
Hence,
\[
 \sum_{n=1}^\infty \lambda^{(1+2\gamma)n}|u_n|^3 \geq
A\|u\|_{\alpha + \gamma}^3,
\]
which concludes the proof.
\end{proof}


\begin{lemma} \label{l:contnnp1}
Let $u(t)$ be a solution to \eqref{model} and $\gamma>0$. Let
$\|u(t)\|^2_\gamma$ be continuous on $[0,\infty)$. Then the function
\begin{equation} \label{e:functionpr}
\sum_{n=1}^\infty \lambda^{2\gamma n} (u_nu_{n+1})(t)
\end{equation}
is continuous on $[0,\infty)$.
\end{lemma}
\begin{proof}
Let $v_n(t)=\lambda^{2\gamma n}(u_nu_{n+1})(t)$. First, note that due to
the Cauchy-Schwarz inequality, the function \eqref{e:functionpr} is less than
or equal to $\|u(t)\|^2_\gamma$ and, consequently, is bounded on every
interval $[a,b]$, $0 \leq a < b$. Let $t_0>0$. Since $\|u(t)\|^2_{\gamma}$ is
continuous at $t=t_0$, it follows that
\[
\lim_{N\to \infty} \limsup_{t\to t_0} \sum_{n=N}^\infty \lambda^{2\gamma n} u_n^2(t)
=0.
\]
Therefore,
\[
\lim_{N\to \infty} \limsup_{t\to t_0} \sum_{n=N}^\infty v_n(t)
=0,
\]
which means that \eqref{e:functionpr} is continuous at $t=t_0$. Indeed, since
$v_n(t)$ is continuous for every $n$, we have
\begin{multline*}
\limsup_{t\to t_0} \left|\sum_{n=1}^\infty v_n(t) - \sum_{n=1}^\infty v_n(t_0)\right|\\
= \lim_{N\to \infty}\limsup_{t\to t_0}  \left|\sum_{n=1}^{N-1}v_n(t)
-\sum_{n=1}^{N-1}v_n(t_0) + \sum_{n=N}^{\infty}v_n(t) -
\sum_{n=N}^{\infty} v_n(t_0) \right|=0.
\end{multline*}
Similarly, the continuity of \eqref{e:functionpr} from the right holds at $t=0$.
\end{proof}

Now we proceed to our main result.

\begin{theorem} \label{thm:blowup}
Let $u(t)$ be a solution to \eqref{model} with $u_n(0) \geq 0$ 
and $\alpha< 1/3$. Then for every $\gamma >0$, there exists
a constant $M(\gamma)$, such that $\|u(t)\|^3_{1/3+\gamma}$ is not locally integrable on $[0,\infty)$, provided $\|u(0)\|_\gamma > M(\gamma)$.
\end{theorem}
\begin{proof}
Since $\|u\|_{\gamma_1} \leq  \|u\|_{\gamma_2}$ for $\gamma_1\leq \gamma_2$,
it is enough to prove the theorem in the case $0<\gamma < \min\{1/3,1-3\alpha\}$.
Given such $\gamma$, let $u(t)$ be a solution to \eqref{model},
such that $\|u(t)\|^3_{1/3+\gamma}$ is integrable on $[0,T]$  for every $T>0$.
We will show that $\|u(0)\|_{\gamma}$ is bounded from above by a constant
dependent of $\gamma$.

Note that $u_n(t) \geq 0$ for all $n, t>0$ due to Theorem~\ref{t:iquality}.
First, we obtain
\begin{equation} \label{integrq}
\begin{split}
\int_0^{T} \sum_{n=1}^\infty \lambda^{(1+2\gamma)n}u_n^2u_{n+1} \, d\tau &\leq
\int_0^{T} \sum_{n=1}^{\infty} \lambda^{(1+2\gamma)n}(u_n^3 +u_{n+1}^3) \, d\tau\\
&\leq 2\int_0^T \left(\sum_{n=1}^{\infty} \lambda^{\frac{2}{3}(1+2\gamma)n}u_n^2\right)^{3/2} \, d\tau\\
&\leq 2 \int_0^{T} \|u(t)\|_{1/3+\gamma}^3 \, d\tau\\
&< \infty,
\end{split}
\end{equation}
for all $T>0$. Thus, 
\[
\sum_{n=1}^{\infty} \lambda^{(1+2\gamma)n}(u_n^2u_{n+1})(t) \qquad \text{and} \qquad 
\sum_{n=1}^{\infty}\lambda^{(1+2\gamma)n}u_n(t)^3
\]
are locally integrable on $[0,\infty)$. 
In addition, since $\alpha <1/3$, we have
\[
\|u(t)\|_{\alpha+\gamma}^2 \leq \|u(t)\|_{1/3+\gamma}^2,
\]
which implies that $\|u(t)\|_{\alpha+\gamma}^2$ is locally integrable
on $[0,\infty)$.

Now note that if $u_n \leq \frac{1}{2}u_{n+1}$, then $u_nu_{n+1}^2 \leq \frac{1}{2}u_{n+1}^3$.
Otherwise, $u_nu_{n+1}^2 \leq 2 u_n^2 u_{n+1}$. Hence,
\begin{equation} \label{eq:ineq1}
u_nu_{n+1}^2 \leq {\ts \frac{1}{2}}u_{n+1}^3 + 2 u_n^2 u_{n+1}, \qquad n \in \mathbb{N}.
\end{equation}
This also implies that
\begin{equation} \label{eq:ineq2}
\begin{split}
u_nu_{n+1} u_{n+2} &\leq {\ts \frac{1}{2}}u_n^2u_{n+1} +
{\ts \frac{1}{2}}u_{n+1}u_{n+2}^2\\
&\leq {\ts \frac{1}{2}}u_n^2u_{n+1} + {\ts \frac{1}{4}}u_{n+2}^3 + u_{n+1}^2u_{n+2},
\end{split}
\end{equation}
for all $n \in \mathbb{N}$.

From \eqref{model} we obtain
\[
\begin{split}
\ddt(u_nu_{n+1}) = &-\nu (\lambda^{2\alpha n} + \lambda^{2\alpha (n+1)})
u_nu_{n+1}\\
& + \lambda^n u_{n-1}^2u_{n+1} - \lambda^{n+1} u_nu_{n+1}^2\\
&  + \lambda^{n+1} u_n^3- \lambda^{n+2} u_nu_{n+1}u_{n+2}\\
& +g_nu_{n+1} +g_{n+1}u_n.
\end{split}
\]
This, together with inequalities \eqref{eq:ineq1} and \eqref{eq:ineq2}, implies
that
\begin{multline*}
\ddt(u_nu_{n+1}) +\nu(1+\lambda^{2\alpha})\lambda^{2\alpha n} u_nu_{n+1}\\
+2\lambda^{n+1}u_n^2u_{n+1} + {\ts \frac{1}{2}}\lambda^{n+2}u_n^2u_{n+1}
+\lambda^{n+2}u_{n+1}^2u_{n+2}\\
\geq \lambda^{n+1}u_n^3 - {\ts \frac{1}{2}}\lambda^{n+1}u_{n+1}^3 -
{\ts \frac{1}{4}} \lambda^{n+2}u_{n+2}^3.
\end{multline*}
Multiplying it by $\lambda^{2\gamma n}$, taking a sum from $1$ to $\infty$,
and integrating between $0$ and $t$, we get 
\begin{equation} \label{eq:long}
\begin{split}
\sum_{n=1}^\infty \lambda^{2\gamma n}&(u_nu_{n+1})(t)-\sum_{n=1}^\infty
\lambda^{2\gamma n}(u_nu_{n+1})(0)\\
&\quad + \nu (1+ \lambda^{2\alpha})\int_{0}^t \sum_{n=1}^\infty \lambda^{2(\alpha+\gamma)n}  u_nu_{n+1} \, d\tau\\
&\quad + (2\lambda  + {\ts \frac{1}{2}}\lambda^2 + \lambda^{1-2\gamma})\int_{0}^t  \sum_{n=1}^\infty \lambda^{(1+2\gamma)n}u_n^2u_{n+1} \, d\tau
\\
&\geq (\lambda-{\ts \frac{1}{2}}\lambda^{-2\gamma}-{\ts \frac{1}{4}}\lambda^{-4\gamma})\int_{0}^t\sum_{n=1}^\infty \lambda^{(1+2\gamma) n} u_n^3 \, d\tau \\
&\geq \frac{\lambda}{4} \int_{0}^t\sum_{n=1}^\infty
\lambda^{(1+2\gamma) n} u_n^3 \, d\tau,
\end{split}
\end{equation}
for all $t\geq 0$.
On the other hand, we have the following equality for the nonlinear term:
\[
\begin{split}
-(B(u,u),A^{\gamma/\alpha}u) &= \sum_{n=1}^{\infty} \lambda^{2\gamma n + n}
u_{n-1}^2 u_n - \sum_{n=1}^{\infty}
 \lambda^{2\gamma n + n+1}u_n^2 u_{n+1} \\
&=c_1 \sum_{n=1}^{\infty} \lambda^{(1+2\gamma)n} u_{n}^2 u_{n+1},
\end{split}
\]
where $c_1= \lambda^{2\gamma+1} - \lambda>0$.
Now, multiplying \eqref{model} by $\lambda^{2\gamma n}u_n$,
taking a sum from $1$ to $\infty$, and integrating between $0$ and $t$,
we obtain

\begin{multline} \label{eq:short}
\|u(t)\|_\gamma^2  - \|u(0)\|_\gamma^2
+2\nu \int_{0}^t \|u(\tau)\|_{\alpha + \gamma}^2 \, d\tau\\
= 2c_1 \int_{0}^t \sum_{n=1}^{\infty} \lambda^{(1+2\gamma)n}u_n^2 u_{n+1}
\, d\tau + 2\int_{0}^t\sum_{n=1}^\infty \lambda^{2\gamma n} g_n u_n \, d\tau,
\end{multline}
for $t\geq 0$. Note that the term with the force is integrable because
$\gamma\leq 1/3$. In particular, \eqref{eq:short} yields that
$\|u(t)\|^2_\gamma$ is continuous on $[0, \infty)$. Denote
\[
H(t) :=  \|u(t)\|_\gamma^2+ c_2\sum_{n=1}^\infty
\lambda^{ 2\gamma n}(u_nu_{n+1})(t),
\]
where $c_2=2c_1/(2\lambda  + \lambda^2/2+ \lambda^{1-2\gamma})$.
Thanks to Lemma~\ref{l:contnnp1}, $H(t)$ is continuous on $[0,\infty)$.
We will show that $H(t)$ is a Lyapunov function, i.e.,
$H(t)$ is always increasing. Moreover, we will see that $H(t)$ blows up in finite
time. Indeed, multiplying \eqref{eq:long} by $c_2$ and adding
\eqref{eq:short}, we get
\begin{multline*}
H(t)-H(0) \geq -2\nu \int_{0}^t \|u(\tau)\|_{\alpha + \gamma}^2 \, d\tau
- \nu c_3
\int_{0}^t \sum_{n=1}^\infty \lambda^{2(\alpha+\gamma)n}  u_nu_{n+1} \, d\tau \\
 + \frac{\lambda c_2}{4}\int_{0}^{t} \sum_{n=1}^\infty   \lambda^{(1+2\gamma) n}u_n^3 \, d\tau,
\end{multline*}
where $c_3=(1+\lambda^{2\alpha})c_2$.
Due to Lemma~\ref{l:jensen}, there exists a constant $A>0$,
such that
\[
 \sum_{n=1}^\infty \lambda^{(1+2\gamma)n}u_n^3 \, \geq
A \|u\|_{\alpha + \gamma}^3.
\]
In addition, the Cauchy-Schwarz inequality implies
\[
\sum_{n=1}^\infty \lambda^{2(\alpha+\gamma)n}  u_nu_{n+1}
\leq \|u\|^2_{\alpha+\gamma}.
\]
Therefore, we obtain
\begin{equation} \label{eq:almlast}
H(t)-H(0) \geq -\nu(2+ c_3)
\int_{0}^t \|u(\tau)\|_{\alpha + \gamma}^2 \, d\tau
+ \frac{A\lambda c_2}{4}\int_{0}^t \|u(\tau)\|_{\alpha + \gamma}^3 \, d\tau,
\end{equation}
for $t\geq 0$.
Note that
\[
\|u(t)\|^2_\gamma  \leq  H(t) \leq (1+c_2) \|u(t)\|^2_\gamma.
\]
In particular,
\[
\|u(t)\|_{\alpha + \gamma} \geq \sqrt{\frac{H(t)}{1+c_2}}. 
\]
Let
\[
M(\gamma) := \frac{8\nu(2+c_3)\sqrt{1+c_2}}{A\lambda c_2}.
\]
Assume that $H(\tau) \geq M^2$ on $[0, t]$ for some $t>0$. Then we have that
$\|u(\tau)\|_{\alpha+\gamma}\geq 8\nu(2+c_3) / (A\lambda c_2) $ on $[0,t]$ and, consequently, \eqref{eq:almlast} yields
\begin{equation} \label{e:ricc}
H(t)-H(0)
\geq c\int_{0}^t H(\tau)^{3/2} \, d\tau,
\end{equation}
where $c=A\lambda c_2/8$.
Now assume that $\|u(0)\|_\gamma > M$. Then $H(0) > M^2$ and
\eqref{e:ricc} holds for some small time $t>0$. Then $\eqref{e:ricc}$
automatically holds for every $t>0$.

Note that \eqref{e:ricc} is a Riccati-type inequality. It is easy to see that $H(t)$
blows up in finite time. Indeed, let $y(t)$ be the solution to the Riccati
equation
\[
y'(t) = c y(t)^{3/2}, \qquad y(0) = {\ts \frac{1}{2}}H(0).
\]
Then for some $t^*>0$, we have that $y(t) \to \infty$ as $t \to t^* -$.
Consider
\[
w(t) = H(t)-y(t).
\]
It is easy to check that the function $w(t)$ satisfies the following integral
inequality:
\[
w(t) -w(0) \geq c \int_{0}^t w(\tau)^{3/2} \, d\tau,
\]
for all $t> 0$, such that $w(\tau)\geq 0$ on $[0,t]$. 
Note that $w(0) >0$ and $w(t)$ is continuous. Thus, $w(t) \geq 0$ for all
$t \in [0,t^*)$.

Now, since $y(t)$ blows up in finite time, $H(t)$ also blows up in finite time, which contradicts the fact that $H(t)$ is continuous on $[0,\infty)$. Hence,
$\|u(0)\|_{\gamma} \leq M$.
\end{proof}

\section{Global attractor} \label{s:attractor}

In Section~\ref{s:apriori} we showed that the dyadic model possesses an
absorbing ball with a radius $R$ larger than $|g|/\nu$.
Let $X$ be a closed absorbing ball.
\[
X:= \{u\in H: |u| \leq R\},
\]
which is compact in $\dw$-metric.
Then for any bounded set $K \subset H$, there exists a time $t_0$, such that
\[
u(t) \in X, \qquad \forall t\geq t_0,
\]
for every Leray-Hopf solution $u(t)$ to \eqref{model} with the initial
data $u(0) \in K$.

We recall the definition of an evolutionary system $\Dc$ from  \cite{CF} (see also \cite{C}). Let
\[
\mathcal{T} := \{ I: \ I=[T,\infty) \subset \mathbb{R}, \mbox{ or } 
I=(-\infty, \infty) \},
\]
and for each $I \subset \mathcal{T}$ let $\mathcal{F}(I)$ denote
the set of all $X$-valued functions on $I$.
A map $\Dc$ that associates to each $I\in \mathcal{T}$ a subset
$\Dc(I) \subset \mathcal{F}$ will be called an evolutionary system if
the following conditions are satisfied:
\begin{enumerate}
\item $\Dc([0,\infty)) \ne \emptyset$.
\item
$\Dc(I+s)=\{u(\cdot): \ u(\cdot -s) \in \Dc(I) \}$ for
all $s \in \mathbb{R}$.
\item $\{u(\cdot)|_{I_2} : u(\cdot) \in \Dc(I_1)\}
\subset \Dc(I_2)$ for all
pairs of $I_1,I_2 \in \Omega$, such that $I_2 \subset I_1$.
\item
$\Dc((-\infty , \infty)) = \{u(\cdot) : \ u(\cdot)|_{[T,\infty)}
\in \Dc([T, \infty)) \ \forall T \in \mathbb{R} \}.$
\end{enumerate}
Let
\[
\begin{split}
R(t)A &:= \{u(t): u(0)\in A, u \in \Dc([0,\infty))\},\\
\Rc(t) A &:= \{u(t): u(0) \in A, u \in \Dc((-\infty,\infty))\}, \qquad A \subset X,
\ t\geq 0.
\end{split}
\]

For $A \subset X$ and $r>0$, denote
$B_{\bullet}(A,r) = \{u: \ d_{\bullet}(A,u) < r\},$ where
$\bullet = \mathrm{s,w}$.
Now we define an attracting set and a global attractor as follows.

\begin{definition}
A set $A \subset X$ is a $\mathrm{d}_{\bullet}$-attracting set
($\bullet = \mathrm{s,w}$) if it uniformly
attracts $X$ in $\mathrm{d}_{\bullet}$-metric, i.e.,
for any $\epsilon>0$ there exists
$t_0$, such that
\[
R(t)X \subset B_{\bullet}(A, \epsilon), \qquad \forall t \geq t_0.
\]
A set $A \subset X$ is invariant if $\Rc(t) A = A$ for all $t\geq 0$.
A set $\mathcal{A}_{\bullet}\subset X$ is a
$\mathrm{d}_{\bullet}$-global attractor if
$\mathcal{A}_{\bullet}$ is a minimal $\mathrm{d}_{\bullet}$-closed
$\mathrm{d}_{\bullet}$-attracting  set.
\end{definition}
The following result was proved in \cite{CF}:
\begin{theorem} \label{c:weakA}
The evolutionary system $\Dc$ always possesses a weak global attractor
$\Aw$.
In addition, if $\Dc([0,\infty))$ is compact in
$C([0, \infty);\Hw)$, then
\begin{enumerate}
\item[(a)]
$ \Aw =\{ u^0: \ u^0=u(0) \mbox{ for some }
u \in \Dc((-\infty, \infty))\}.$
\item[(b)] $\Aw$ is the maximal invariant set.
\end{enumerate}
\end{theorem}

For the dyadic model, we define $\Dc$ in the following way.
\[
\begin{split}
\Dc([T,\infty)) := &\{u: u(\cdot)
\mbox{ is a Leray-Hopf solution on } [T,\infty)\\
& \ \mbox{ and } u(t) \in X \ \forall t \in [T,\infty)\},
\qquad T \in \mathbb{R},
\end{split}
\]
\[
\begin{split}
\Dc((\infty,\infty)) := & \{u: u(\cdot)
\mbox{ is a Leray-Hopf solution on } (-\infty,\infty)\\
& \ \mbox{ and } u(t) \in X \ \forall t \in (-\infty,\infty)\},
\end{split}
\]
where $X$ is the phase space defined in the beginning of the section.
Clearly, $\Dc$ satisfies properties (1)--(4). Then Theorem~\ref{c:weakA} immediately 
yields that the weak global attractor $\Aw$ exists. In order to infer
that $\Aw$ is the maximal invariant set, we need the following result.

\begin{lemma} \label{l:compact}
$\Dc([0,\infty))$ is compact in $C([0, \infty);\Hw)$.
\end{lemma}
\begin{proof}
Take any sequence $u^k \in \Dc([0,\infty))$. 
First, note that
\[
u^{k}_n(t) \leq R, \qquad \forall n, k, t\geq 0.
\]
Therefore,
\[
|u^k_n(t) - u^k_n(s)|
\leq (\nu\lambda^{2\alpha n}R +\lambda^n R^2 + \lambda^{n+1}R^2 + g_n)
|t-s|,
\]
for all $n$, $k$, $t\geq 0$, $s\geq0$. Thus,
\[
\dw(u^{k}(t), u^{k}(s)) = \sum_{n=1}^\infty \frac{1}{2^{(n^2)}}
\frac{|u^k_n(t)-u^k_n(s)|}{1 + |u^k_n(t)-u^k_n(s)|} \leq c|t-s|,
\]
for some constant $c$ independent of $k$. Hence, $\{u^{k}\}$ is an
equicontinuous sequence
of functions in $C([0,\infty);\Hw)$ with bounded initial data. Therefore,
Ascoli-Arzela theorem implies that
$\{u^k\}$ is relatively compact in $C([0,T]; \Hw)$ for all $T>0$. Using
a diagonalization process, we obtain that $\{u^k\}$ is relatively compact in
$C([0,\infty); \Hw)$. Hence, there exists a weakly continuous $H$-valued
function $u(t)$ on $[0,\infty)$, such that
\begin{equation} \label{weakconv}
u^{k_j} \to u \qquad \mbox{as} \qquad k_j \to \infty  \qquad \mbox{in}
\qquad  C([0,\infty); \Hw),
\end{equation}
for some subsequence $k_j$.
In particular,
\[
|u(t)| \leq \liminf_{k_j \to \infty} |u^{k_j}(t)| \leq R, \qquad t \geq 0,
\]
i.e., $u(t) \in X$ for all $t \geq 0$.

In addition, since $u^{k_j}(t)$ is a solution to \eqref{model}, we have
\[
u^{k_j}_n(t) = u^{k_j}_n(0) + \int_0^t (-\nu \lambda^{2\alpha n}u^{k_j}_n +
\lambda^n (u^{k_j}_{n-1})^2 - \lambda^{n+1} u^{k_j}_n u^{k_j}_{n+1} +g_n) \, d\tau,
\]
for all $n$. Taking the limit as $k_j \to \infty$, we obtain
\[
u_n(t) = u_n(0) + \int_0^t (-\nu \lambda^{2\alpha n}u_n +
\lambda^n u_{n-1}^2 - \lambda^{n+1} u_n u_{n+1} +g_n) \, d\tau,
\]
for all $n$.
Since $u_n(t)$ is continuous, $u_n \in C^1([0,\infty))$ and satisfies
\eqref{model}.

In order to infer that $u \in \Dc([0,\infty))$,
it remains to prove that $u(t)$ satisfies the energy inequality.
Note that  $|u^{k_j}(t)| \to |u(t)|$ as $k_j \to \infty$ a.e. in $[0,\infty)$.
Since $u^{k_j} \in \Dc([0,\infty))$, it satisfies the energy inequality
starting from any $t_0$ that is not in the exceptional set of measure zero.
Let $Ex$ be the union of the exceptional sets for all $u^{k_j}$. Note that $Ex$
is of measure zero.
Take any $t_0 \notin Ex$ for which $|u^{k_j}(t_0)| \to |u(t_0)|$ as $k_j \to \infty$.
Then 
\[
|u^{k_j}(t)|^2 + 2\nu \int_{t_0}^t \|u^{k_j}(\tau)\|^2 \, d\tau \leq
|u^{k_j}(t_0)|^2 + 2\int_{t_0}^t (g, u^{k_j}(\tau)) \, d\tau,
\]
for all $t \geq t_0$. Hence,
\[
|u^{k_j}(t)|^2 + 2\nu \int_{t_0}^t \sum_{n\leq N} \lambda^{2\alpha n} u^{k_j}_n(\tau)^2 \, d\tau \leq
|u^{k_j}(t_0)|^2 + 2\int_{t_0}^t (g, u^{k_j}(\tau)) \, d\tau.
\]
Since $u^{k_j}(t) \to u(t)$ weakly in $H$ as $k_j \to \infty$
for all time $t\geq 0$, we have that
\[
|u(t)|^2 + 2\nu \int_{t_0}^t \sum_{n\leq N} \lambda^{2\alpha n} u_n(\tau)^2 \, d\tau \leq
|u(t_0)|^2 + 2\int_{t_0}^t (g, u(\tau)) \, d\tau.
\]
Finally, taking the limit as $N \to \infty$ and using Levi's convergence theorem, we obtain
\[
|u(t)|^2 + 2\nu \int_{t_0}^t \|u(\tau)\|^2 \, d\tau \leq
|u(t_0)|^2 + 2\int_{t_0}^t (g, u(\tau)) \, d\tau,
\]
for all $0 \leq t_0 \leq t$, $t_0$ a.e. in $[0,\infty)$. Hence, $u \in \Dc([0,\infty))$, which concludes the proof.
\end{proof}

Now Theorem~\ref{c:weakA} implies that the weak global attractor
$\Aw$ is the maximal invariant set that consists of the points that belong to
complete trajectories, i.e., trajectories in $\Dc((\infty,\infty))$. Moreover, 
using \eqref{e:gensol}, one
can show that $u_n\geq 0$ for every $u\in \Aw$. Consider now the case
$\alpha<1/3$. It is easy to show that for every
$\gamma \in (0, 1-3\alpha)$, we can take $g_1$ large enough, so that
for every solution $u(t)$ and every $t \geq 0$, we have $|u(\tau)|>M(\gamma)$
for some $\tau \in [t,t+1]$.
Thanks to Theorem~\ref{thm:blowup}, this means that
$\Aw$ is not bounded in $H^{1/3+\gamma}$.
It is an open question whether $\Aw$ is bounded in $V$.

We will now proceed to study the question whether $\Aw$ is also
a strong global attractor.

\begin{theorem}
Let $\alpha \geq 1/2$. Then every Leray-Hopf solution $u(t)$ of \eqref{model}
satisfies the energy equality
\[
|u(t)|^2 + 2\nu \int_{t_0}^t \|u(\tau)\|^2 \, d\tau =
|u(t_0)|^2 + 2\int_{t_0}^t (g, u(\tau)) \, d\tau, \qquad 0\leq t_0 \leq t.
\] 
\end{theorem}
\begin{proof}
Let $u(t)$ be a Leray-Hopf solution of \eqref{model}.
Thanks to the energy inequality \eqref{ee},
$\|u(t)\|^2$ is locally integrable.
Then we obtain
\begin{equation*}
\begin{split}
\int_{t_0}^{t} \sum_{n=1}^\infty \lambda^{n}u_n^2u_{n+1} \, d\tau &\leq
\int_{t_0}^{t} \sum_{n=1}^\infty \lambda^{n}(u_n^3 +u_{n+1}^3) \, d\tau\\
&\leq 2\sup_{s\in[t_0,t]}|u(s)|\int_{t_0}^{t}\sum_{n=1}^\infty \lambda^n u_n^2 \, d\tau\\
&\leq 2R\int_{t_0}^{t}\|u(\tau)\|^2 \, \, d\tau\\
&< \infty,
\end{split}
\end{equation*}
for $0 \leq t_0 \leq t$. Hence,
\[
\int_{t_0}^t \lambda^{n+1} u_n^2u_{n+1} \, d\tau \to 0, \qquad \text{as}
\qquad  n \to \infty.
\]
Multiplying \eqref{model} by $u_n$,
taking a sum from $1$ to $N$, and integrating between $t_0$ and $t$,
we obtain
\begin{multline*}
\sum_{n=1}^N u_n(t)^2  - \sum_{n=1}^N u_n(t_0)^2
 +2\nu \int_{t_0}^t \sum_{n=1}^N \lambda^{2\alpha n}u_n(\tau)^2
 \, d\tau \\
= - 2\int_{t_0}^t  \lambda^{N+1} u_N^2
u_{N+1} \, d\tau  +2\int_{t_0}^t\sum_{n=1}^N  g_n u_n.
\end{multline*}
Finally, taking the limit as $N\to \infty$, we arrive at
\[
|u(t)|^2 + 2\nu \int_{t_0}^t \|u(\tau)\|^2 \, d\tau =
|u(t_0)|^2 + 2\int_{t_0}^t (g, u(\tau)) \, d\tau, \qquad 0\leq t_0 \leq t.
\]
\end{proof}

In \cite{CF} it was proved that the asymptotic compactness of the
dynamical system $\Dc$ implies that the strong global attractor $\As$
exists, is strongly compact, and coincides with $\Aw$. In the case
where the evolutionary system consists of the Leray-Hopf weak solutions
to the 3D NSE, the continuity of the complete trajectories, i.e. trajectories on
$\Aw$, implies the asymptotic
compactness of $\Dc$ (see also \cite{B1} and  \cite{R} for similar results).
In \cite{C} this result was proved for an abstract evolutionary system
satisfying the energy inequality. It immediately implies the following.

\begin{corollary} \label{col:attr}
Let $\alpha \geq 1/2$.  Then $\Aw$
is a strongly compact strong global attractor.
\end{corollary} 

Note that if $\alpha \geq 1/2$, then, thanks to Theorem~\ref{t:gregularity},
for every initial datum in $H$ there exists a
regular solution on $[0,\infty)$. Moreover, it can be shown that such a solution
is unique in the class of all Leray-Hopf solutions.
Hence, Corollary~\ref{col:attr} can also
be obtained using a classical theory of semiflows. It is an open question
whether the continuity of the complete trajectories and, consequently,
the existence of the strong compact global attractor holds for $\alpha <1/2$.

\section*{acknowledgment}
The author would like to thank Susan Friedlander and anonymous referee
for helpful comments and suggestions.

\end{document}